\documentclass[11pt]{article}
\usepackage{amsxtra,amssymb,amsmath,enumerate,graphicx,bbm,amsthm,color}
\usepackage[active]{srcltx}
\usepackage{t1enc}
\usepackage{aurical}
 \usepackage{ulem}
 
 \usepackage{geometry}
\def\be{\begin{eqnarray}}
\def\ee{\end{eqnarray}}
\def\b*{\begin{eqnarray*}}
\def\e*{\end{eqnarray*}}
\def\Hb{\mathbf{H}}
\def\Sb{\mathbf{S}}
\def\Lb{\mathbf{L}}
\def\Ab{\mathbf{A}}
\def\Xb{\overline{X}}

\def\Tc{\mathcal{T}}
\def\Fc{\mathcal{F}}
\def\Ec{\mathcal{E}}
\def\Mc{\mathcal{M}}
\def\Oc{\mathcal{O}}
\def\Uc{\mathcal{U}}
\def\Qc{\mathcal{Q}}
\def\eps{\varepsilon}

\def\E{\mathbb{E}}
\def\R{\mathbb{R}}
\def\P{\mathbb{P}}
\def\Q{\mathbb{Q}}
\def\F{\mathbb{F}}
\def\N{\mathbb{N}}
\def\vs#1{\vspace{#1mm}}
\def\esssup{{\rm ess\!}\sup}
\def\x{\times}
\def \proof{{\noindent \bf Proof. }}
\def \ep{\hbox{ }\hfill$\Box$}
 \def\reff#1{{\rm(\ref{#1})}}
 \def \1{{\bf 1}}
 \def\as{{\rm a.s.}}
 \def\Ecg{\Ec^{g}}
\def\Tcd{\Tc_{2}}

\def\HTC{{\rm Tc}}
\def\HS{{\rm S}}
\def\HSLD{{\rm Sld}}
\def\hTC{{\rm (Tc)}}
\def\hS{{\rm (S)}}
\def\hSLD{{\rm (Sld)}}

\newtheorem{Theorem}{Theorem}[section]

\newtheorem{Proposition}{Proposition}[section]
\newtheorem*{Assumption}{Assumption}
\newtheorem{Lemma}{Lemma}[section]
\newtheorem{Corollary}{Corollary}[section]
\newtheorem{Remark}{Remark}[section]

\begin{document}

\title{A general Doob-Meyer-Mertens decomposition for $g$-supermartingale systems\thanks{We are very grateful to N.~El Karoui for discussions we had on the first version of this paper.}}

\author{Bruno Bouchard\footnote{CEREMADE and CREST-LFA, University Paris-Dauphine and ENSAE-ParisTech. ANR Liquirisk. bouchard@ceremade.dauphine.fr} 
\and Dylan Possama\"i\footnote{CEREMADE,  University Paris-Dauphine. possamai@ceremade.dauphine.fr} 
\and Xiaolu Tan\footnote{CEREMADE,  University Paris-Dauphine. tan@ceremade.dauphine.fr, the author gratefully acknowledges the financial support of the ERC 321111 Rofirm, the ANR Isotace, and the Chairs Financial Risks (Risk Foundation, sponsored by Soci\'et\'e G\'en\'erale) and Finance and Sustainable Development (IEF sponsored by EDF and CA).}}

\date{First version: May 2015. This version: {July} 2015}

\maketitle

\begin{abstract}
	We provide a general Doob-Meyer decomposition for $g$-supermartingale systems,
	which does not require any right-continuity on the system, {nor that the filtration is quasi left-continuous}.
	In particular, it generalizes the Doob-Meyer decomposition of Mertens \cite{mertens1972theorie} for classical supermartingales,
	as well as Peng's \cite{peng1999monotonic} version for right-continuous $g$-supermartingales.
	As examples of application, 
	we prove an optional decomposition theorem for $g$-supermartingale systems,
	and also obtain a general version of the well-known dual formation for BSDEs with constraint on the gains-process, using very simple arguments.

\end{abstract}

\vspace{5mm}


\vspace{5mm}

\noindent {\textbf{MSC Classification (2010)}: 60H99.

 \section{Introduction}

The Doob-Meyer decomposition is one of the fundamental result of the general theory of processes, in particular when applied to the theory of optimal control, see El Karoui \cite{stflour}. Recently, it has been pointed out by Peng \cite{peng1999monotonic} that it also holds in the semi-linear context of the so-called $g$-expectations. 	
	Namely, let $(\Omega, \Fc, \P)$ be a probability space equipped with a $d$-dimensional Brownian motion $W$,
	as well as the Brownian filtration $\F = (\Fc_t)_{t \ge 0}$,
	let $g: (t,\omega, y, z) \in \R^+ \x \Omega \x \R \x \R^d \longrightarrow \R$  be some function, progressively measurable in $(t,\omega)$ and Lipschitz in $(y,z)$,
	and $\xi \in \Lb^{2}(\Fc_{\tau})$ for some stopping time $\tau$. We define $\Ec^{g}_{\cdot,\tau}[\xi]:=Y_{\cdot}$ 
	where $(Y,Z)$ solves the backward stochastic differential equation 
	$$
		-dY_{t} = g_{t}(Y_{t},Z_{t})dt - Z_{t}\cdot dW_{t},
		\  \mbox{on}\ 
		[0,\tau],
	$$  
	with terminal condition $Y_{\tau}=\xi$.  Then, an optional process $X$ is said to be a (strong) $\Ec^{g}$-supermartingale if for all stopping times $\sigma\le \tau$ we have $X_{\tau}\in \Lb^{2}(\Fc_{\tau})$ and  $X_{\sigma}\ge \Ec^{g}_{\sigma,\tau}[X_{\tau}]$ almost surely. 
	When $X$ is right-continuous, it admits a unique decomposition of the form
	$$
		-dX_{t} = g_{t}(X_{t},Z^{X}_{t})d t - Z^{X}_{t}\cdot dW_{t} + dA^{X}_{t},$$ 
	in which $Z^{X}$ is a square integrable and predictable process, 
	and $A^{X}$ is   non-decreasing predictable. 
	See \cite{peng1999monotonic} and \cite{coquet2002filtration,hu2008representation,ma2010quadratic}. 
	In particular, when $g\equiv 0$, this is the classical Doob-Meyer decomposition in a Brownian filtration framework. 

\vspace{0.5em}
\noindent As fundamental as its classical version, this result was used by many authors in various contexts : backward stochastic differential equation with constraints  \cite{bouchard2014regularity,klimsiak,peng2010reflected}, minimal supersolutions under non-classical conditions on the driver \cite{drapeau2013minimal,heyne2014minimal}, minimal supersolutions under volatility uncertainty \cite{cheng2013optimal,samuel2013minimal,mpp,mpz,pos1,possamai2013second,sontouzha11,SoToZh09a}, backward stochastic differential equations with weak terminal conditions \cite{reveillac2012bsdes}, etc.

\vspace{0.5em}
\noindent However, it is limited to right-continuous $\Ec^{g}$-supermartingales,
	while the right-continuity might be very difficult to prove, if even correct. 
	The method generally used by the authors is then to work with the right-limit process, which is automatically right-continuous, but they then face important difficulties in trying to prove that it still shares the dynamic programming principle of the original process.
	This was sometimes overcome to the price of stringent assumptions,
	which are often too restrictive, in particular in the context of singular optimal control problems.

\vspace{0.5em}
\noindent In the classical case, $g \equiv 0$, it is well known that we can avoid these technical difficulties by appealing to the version of the Doob-Meyer decomposition for supermartingales with only right and left limits, see El Karoui \cite{stflour}.
	It has been established by Mertens \cite{mertens1972theorie}, Dellacherie and Meyer \cite[Vol.~II, Appendice~1]{DellacherieMeyer} provides an alternative proof.
	Unfortunately, such a result has not been available so far in the semi-linear context.

 \vspace{0.5em}
\noindent This paper fills this gap\footnote{After completing this manuscript, we discovered \cite{grigo} that was issued at the same time. In this paper, the authors prove the existence of reflected BSDEs for barriers with only right-limits, from which they can infer a similar Doob-Meyer-Mertens decomposition as the one proved in the current paper. Their decomposition is less general, in terms of integrability conditions and assumptions on the filtration.  On the other hand, we do not provide any comparable existence result for reflected BSDEs with only right-limited barriers {(see however our companion paper \cite{BPTZ}, where a general existence result for reflected BSDEs with c\`adl\`ag obstacles is given)}. Also, our technic of proof is quite different.} and  provides a version  {\sl à la} Mertens of the Doob-Meyer decomposition of $\Ec^{g}$-supermartingales. 
	By following the arguments of Mertens \cite{mertens1972theorie}, we first show that a supermartingale associated to a general family of semi-linear (non-expansive) and time consistent expectation operators can be corrected into a right-continuous one by subtracting the sum of the previous jumps on the right. 
	Applying this result to the $g$-expectation context, together with the  decomposition of \cite{peng1999monotonic}, we then obtain a decomposition for the original $\Ecg$-supermartingale, even when it is not right-continuous. The same arguments apply to $g$-expectations defined  on $\Lb^{p}$, $p>1$, and more general filtrations than the Brownian one considered in \cite{peng1999monotonic}, {in particular we shall not assume that the filtration is quasi left-continuous}. This is our Theorem \ref{thm: doob-mertens decompo bsde} below.  The only additional difficulty is that the decomposition for right-continuous processes has to be extended first. This is done by using the fact that it can naturally be obtained by considering the BSDE reflected from below on the $\Ec^{g}$-supermartingale and by using recent technical extensions of the seminal paper El Karoui et al.~\cite{KKPPQ}, see Proposition \ref{prop:doob-general} below. Then, using classical results of the general theory of stochastic processes, we can even replace the notion of supermartingale by that of supermartingale systems, for which an optional aggregation process can be easily found, see El Karoui \cite{stflour} for the classical case $g\equiv 0$.

  \vspace{0.5em}
\noindent These key statements aim not only at extending already known results to much more general contexts, 
	but also at simplifying many difficult arguments recently encountered in the literature. 
	We provide two illustrative examples. First, we prove a general optional decomposition theorem for  $g$-supermartingales. To the best of our knowledge, such a decomposition was not obtained before. Then, we show how a general duality  for the minimal super-solution of a backward stochastic differential equation with constraint on the gains-process can be obtained. This is an hold problem, but we obtain it  in a framework that could not be considered in the literature before, compare with \cite{bouchard2014regularity,cvitanic1998backward}. In both cases, these a-priori difficult results turn out to be easy consequences of our main Theorem \ref{thm: doob-mertens decompo bsde}, whenever right continuity {\it per se} is irrelevant.
	
	\vspace{3mm}
 
	\noindent {\bf Notations:} 
	$\mathrm{(i)}$ In this paper, 
	$(\Omega,\Fc,\P)$ is a complete probability space, endowed with a filtration $\F=(\Fc_{t})_{t\ge 0}$ satisfying the usual conditions. {Note that we do not assume that the filtration is quasi left-continuous.}
	
	\vspace{1mm}
	
	\noindent $\mathrm{(ii)}$
	We fix a {\sout{fixed}} time horizon $T>0$ throughout the paper, 
	and denote by $\Tc$ the set of stopping times $\as$~less than $T$. 
	We shall also make use of the set  $\Tc_{\sigma}$ of stopping times $\tau\in\Tc$ $\as$~greater than $\sigma\in \Tc$. For ease of notations, let us say that  $(\sigma,\tau)\in \Tcd$ if $\sigma\in \Tc$ and $\tau \in \Tc_{\sigma}$. 
	
	\vspace{1mm}
	
	\noindent $\mathrm{(iii)}$
	Let $\sigma \in \Tc$, conditional expectations or probabilities given $\Fc_{\sigma}$ are simply denoted by $\E_{\sigma}$ and $\P_{\sigma}$. Inequalities between random variable are taken in the $\as$~sense unless something else is specified. If $\mathbb Q$ is another probability measure on $(\Omega,\Fc)$, which is equivalent to $\P$, we will write $\Q\sim\P$. 
	
	\vspace{1mm}
	
	\noindent $\mathrm{(iv)}$ 
	For any sub-$\sigma$-field $\mathcal G$ of $\mathcal F$, $\Lb^0(\mathcal G)$ denotes the set of random variables on $(\Omega,\Fc)$ which are in addition $\mathcal G$-measurable. Similarly, for any $p\in (0,\infty]$, and any probability measure $\Q$ on $(\Omega, \mathcal F)$, we let $\Lb^{p}(\mathcal G,\Q)$ be the collection of real-valued $\mathcal G$-measurable random variables with absolute value admitting a $p$-moment under $\Q$. For ease of notations, we denote $\Lb^{p}(\mathcal G):=\Lb^p(\mathcal G,\P)$ and also $\Lb^p:=\Lb^p(\Fc)$. These spaces are endowed with their usual norm. 

	\vspace{1mm}

	\noindent $\mathrm{(v)}$
	For $p \in (0, \infty]$, we denote by ${\bf X}^{p}$ (resp.~${\bf X}^{p}_{r}$, ${\bf X}_{ \ell r}^{p}$) the collection of all optional  processes $X$  such that $X_{\tau}$ lies in $\Lb^{p}(\Fc_{\tau})$ for all $\tau \in \Tc$  (resp.~and such that $X$ admits right-limits, and such that $X$ admits right- and left-limits).
	We denote by $\Sb^p$ the  set of all  c\`adl\`ag, $\F$-optional processes $Y$, 
	such that $\sup_{0 \le t \le T} Y_t \in \Lb^p$,
	and by $\Hb^{p}$ the set of all predictable $d$-dimensional processes $Z$ such that
	$$\E\left[\left(\int_0^T|Z_s|^2ds\right)^{\frac p2}\right]<+\infty.$$ 
	Finally, we denote by  $\Ab^{p}$ the set of all non-decreasing predictable processes $A$ such that $A_0 = 0$ and 
	$A_T \in \Lb^{p}$.

	\vspace{1mm}

	\noindent $\mathrm{(vi)}$	For any $d\in\mathbb N\backslash\{0\}$, we will denote by $x\cdot y$ the usual inner product of two elements $(x,y)\in\mathbb R^d\times\mathbb R^d$. We will also abuse notation and let $|x|$ denote the Euclidean norm of any $x\in\R^d$, as well the associated operator norm of any $d\times d$ matrix with real entries.

\section{Stability of $\Ec$-supermatingales under Mertens's re-gularization}

In this section, we provide an abstract regularization result for supermartingales associated to a family of semi-linear non-expansive and time consistent conditional  expectation operators (see below for the exact meaning we give to this, for the moment, vague appellation). It states that we can always modify a supermartingale with right-limits so as to obtain a right-continuous process which is still a supermatingale. This was the starting point of Mertens's proof of the Doob-Meyer decomposition theorem for supermatingales (in the classical sense) with only right-limits. Our proof actually mimics the one of Mertens \cite{mertens1972theorie}. This abstract formulation has the merit to point out the key ingredients that are required for it to go through, in a non-linear context.  It will then be applied to $g$-expectation operators, in the terminology of Peng \cite{peng1997}, to obtain our Doob-Meyer type decomposition, which is the main result of this paper. 

\subsection{Semi-linear time consistent expectation operators}

Let $p\in(1,+\infty]$. Throughout the paper, $q$ will denote the conjugate of $p$
(i.e. $p^{-1} + q^{-1} = 1$). Then, we define a non-linear conditional expectation operator as a family $\Ec=\{\Ec_{\sigma,\tau},\;(\sigma,\tau) \in \Tcd\}$ of maps
$$
\Ec_{\sigma,\tau} : \Lb^{p}(\Fc_{\tau})\longmapsto \Lb^{p}(\Fc_{\sigma}),\; \mbox{ for } (\sigma,\tau) \in \Tcd.
$$
One needs it to satisfy certain structural and regularity properties. Let us start with the notions related to time consistency. 
\begin{Assumption}[\HTC] Fix $(\tau_{i})_{i\le 3}\subset \Tc$ such that $\tau_{1}\vee \tau_{2}\le \tau_{3}$. Then, 
\begin{enumerate}[{\rm (a)}]  
\item $\Ec_{\tau_{1},\tau_{1}}$ is the identity.
\item $\Ec_{\tau_{1},\tau_{2}}\circ\Ec_{\tau_{2},\tau_{3}}=\Ec_{\tau_{1},\tau_{3}}$, if $\tau_{1}\le \tau_{2}$.
\item $\Ec_{\tau_{1},\tau_{3}}[\xi]=\Ec_{\tau_{2},\tau_{3}}[\xi]$ $\as$~on $\{\tau_{1}=\tau_{2}\}$,  for all $\xi \in \Lb^{p}(\Fc_{\tau_{3}})$.
\end{enumerate}
\end{Assumption}
\noindent We also need some regularity with respect to monotone convergence. 
\begin{Assumption}[\HS]   Fix $(\sigma,\tau)\in \Tcd$.
\begin{enumerate}[{\rm (a)}]  
\item Fix $s\in [0,T)$ and   $\xi \in \Lb^{0}(\Fc_s)$. Let $(s_{n})_{n\geq 1}\subset [s,T]$ decrease to $s$ and $(\xi_{n})_{n\ge 1}$ be such that $\xi_{n}\in \Lb^{p}(\Fc_{s_{n}})$ for each $n$, $(\xi_{n}^{-})_{n\ge 1}$ is bounded in $\Lb^{p}$, and $\xi_{n}\longrightarrow \xi$ $\as$ as $n\longrightarrow \infty$, then 
$$
\limsup_{n\to \infty} \Ec_{s,s_{n}}[\xi_{n}]\ge \xi.
$$
\item Let $(\sigma_{n})_{n \ge 1}\subset \Tc$ be a decreasing sequence which converges {\rm a.s.~}to  $\sigma$ and s.t.~$\sigma_n \le \tau$ \as~for all $n \ge 1$. Fix $\xi \in \Lb^{p}(\Fc_{\tau})$. Then, 
$$
\limsup_{n\to \infty} \Ec_{\sigma_{n},\tau}[\xi]\ge  \Ec_{\sigma,\tau}[\xi].
$$
\item Let $(\xi_{n})_{n \ge 1}\subset \Lb^{p}(\Fc_{\tau})$ be a non-decreasing sequence which converges {\rm a.s.}~to $\xi\in \Lb^{p}$.    Then, 
$$
\limsup_{n\to \infty} \Ec_{\sigma,\tau}[\xi_{n}]\ge  \Ec_{\sigma,\tau}[\xi].
$$
\end{enumerate}
\end{Assumption}

\noindent The idea that $\Ec$ should be semi-linear and non-expansive is encoded in the following.
 
 	Let $\Q^1$, $\Q^2$ be two probability measures on $(\Omega, \Fc)$ and $\tau \in \Tc$,
	we define the concatenated probability measure $\Q^1 \otimes_{\tau} \Q^2$ on $(\Omega, \Fc)$ by
	$$
		\E^{\Q^1 \otimes_{\tau} \Q^2} \big[ \xi \big]
		~:=~
		\E^{\Q^1} \big[ \E^{\Q^2} \big[ \xi  \big| \Fc_{\tau} \big] \big],
		~~\mbox{for all bounded measurable variable}~\xi.
	$$
\begin{Assumption}[\HSLD]    
There is a family $\Qc$ of $\P$-equivalent probability measures    such   that: 
\begin{itemize}
\item $\E\left[\left|\frac{d\Q}{d\P}\right|^{q}+\left|\frac{d\Q}{d\P}\right|^{1-q}\right]\le L$ for all $\Q\in \Qc$, for some $L>1$.
\item $\Q^{1} \otimes_{\tau} \Q^{2} \in \Qc$, for all $\Q^{1},\Q^{2}\in \Qc$ and  $\tau \in \Tc$.
\item For all $(\sigma,\tau)\in \Tcd$ and $(\xi,\xi') \in \Lb^{p}(\Fc_{\tau})\x\Lb^{p}(\Fc_{\tau})$ there exists $\Q\in \Qc$  and a $[L^{-1},1]$-valued  $\beta\in \Lb^{0}(\Fc)$ satisfying
$$
\Ec_{\sigma,\tau}[\xi]\le\Ec_{\sigma,\tau}[\xi']+\E^{\Q}_{\sigma}[\beta(\xi-\xi')].
$$ 
\end{itemize}
\end{Assumption}

\noindent Let us comment   this last condition. Assume that $(\Q,\beta)$ is the same for $(\xi,\xi')$ and $(\xi',\xi)$. Then,  inverting the roles of $\xi$ and $\xi'$, it indeed  says that 
$$
\Ec_{\sigma,\tau}[\xi]-\Ec_{\sigma,\tau}[\xi']= \E^{\Q}_{\sigma}[\beta (\xi-\xi')]. 
$$
Otherwise stated, in this case, the operator $\Ec$ can be linearized as each point. However, the linearization, namely $(\Q,\beta)$, depends in general on $(\xi,\xi'),\sigma$ and $\tau$, so that it is not a linear operator. Thus the label semi-linear. 

\vspace{0.5em}
\noindent In any case, it is non-expansive in the sense that $\Ec_{\sigma,\tau}[\xi]-\Ec_{\sigma,\tau}[\xi']\le \E^{\Q}_{\sigma}[|\xi-\xi'|]$, since  $\beta\le 1$. Moreover,    $\Ec_{\sigma,\tau}[\xi]\le \Ec_{\sigma,\tau}[\xi']$ whenever $\xi\le \xi'$ $\as$,  and with strict inequality on $\{\P_{\sigma}[\xi<\xi']>0\}$, since $\beta>0$.

\subsection{Stability by regularization on the right}
 
 Before stating the main result of this section, one needs to define the notion of $\Ec$-supermartingales.
 \vs2
  
\noindent 
We say that $X$ is a $\Ec$-supermatingale if $X\in {\bf X}^{p}$ and  $X_{\sigma}\ge \Ec_{\sigma,\tau}[X_{\tau}]$ $\as$~for all $(\sigma,\tau)\in \Tcd$. We say that it is a local $\Ec$-supermatingale if there exists a non-decreasing sequence of stopping times $(\vartheta_{n})_{n\geq 1}$ s.t. $X_{\sigma\wedge \vartheta_{n}}\ge \Ec_{\sigma\wedge \vartheta_{n},\tau\wedge \vartheta_{n}}[X_{\tau\wedge \vartheta_{n}}]$ for all $(\sigma,\tau)\in \Tcd$ and $n\ge 1$, and $\vartheta_{n}\uparrow \infty$ $\as$~as $n\longrightarrow \infty$.  

\begin{Lemma}\label{lem: X+I supermart abstract} Let Assumptions $\hTC$, $\hS$ and $\hSLD$ hold. Let $X\in {\bf X}^{p}_{r}$ be a $\Ec$-supermartingale such that $(X^{-}_{t})_{t\le T}$ is bounded in $\Lb^{p}$.  Define the process $I$ by 
\begin{equation}\label{eq: def I}
I_{t}:=\sum_{s<t} (X_{s}-X_{s+}),\;t\le T.
\end{equation}
Then, $I$ is non-decreasing, left-continuous and belongs to ${\bf X}^{\frac1p}$. Moreover, $\Xb:=X+I$ is a right-continuous local $\Ec$-supermatingale.
\end{Lemma}

\proof We split the proof in several steps. As already mentioned, we basically only check that the arguments of Mertens \cite{mertens1972theorie} go through under our assumptions.

\vspace{0.5em}
\noindent{(a)\sl~ $\Xb$ is right-continuous}. 
	Indeed, for every $t \in [0,T)$, one has
	$$\Xb_{t+} ~=~ X_{t+} + I_{t+} ~=~ X_{t+}+I_t+(X_t-X_{t+}) ~=~ X_t+I_t.$$

\vspace{0.5em}
\noindent{(b)\sl~ Jumps from the right are non-positive, i.e. $X_t \ge X_{t+}$ for each $t \in [0, T)$, so that $I$ is non-decreasing, and $X_{\sigma+}\ge \Ec_{\sigma,\tau}[X_{\tau}]$ for all $(\sigma,\tau)\in \Tcd$ with $\sigma<\tau$.}  

\vspace{0.5em}
\noindent By the $\Ec$-supermartingale property,   
$X_{t}\ge \Ec_{t,t+h}[X_{t+h}]$ for any $h\in (0, T-t]$ and $t<T$. Since $X_{t+h}\longrightarrow X_{t+}$ as $h\downarrow 0$ and $(X_{t+h}^{-})_{h}$ is bounded in $\Lb^{p}$, 
it follows from \hS(a) that $X_t \ge X_{t+}$.  Similarly, $X_{\sigma+}\ge \Ec_{\sigma,\tau}[X_{\tau}]$ as a consequence of \hS(b).

\vspace{0.5em}
\noindent{(c)\sl~   Let $k \in \N$, $\eps > 0$, and $(\sigma_{i})_{i\le k}\subset \Tc$ be the non-decreasing sequence of stopping times which exhausts the first $k$ jumps from the right of $X$ of size bigger than $\eps$ $($recall that $X$ admits right-limits$)$.
Denote
\be \label{eq:IX_eps_k}
I^{\eps, k}_t :=\sum_{i=1}^{k}  (X_{\sigma_{i}}-X_{\sigma_{i}+}) \1_{ \sigma_{i} < t },
\quad\mbox{and}\quad
\Xb^{\eps, k}_t := X_t + I^{\eps, k}_t,
\ee
then $\Xb^{\eps,k}$ is still a $\Ec$-supermartingale.
}

\vspace{0.5em}
\noindent Note that we can always assume that there is $\as$ at least $k$ jumps, as we can always add jumps of size $0$ at $T$. We shall use the conventions $\sigma_{0}=0$ and $\sigma_{k+1}=T$. 
The proof proceeds by induction and requires several steps. 
For ease of notation, we omit the superscript $(\eps, k)$ in $(\Xb^{\eps, k}, I^{\eps, k})$ and write $(\Xb, I)$ in this part (that is in item (c) only).

\vspace{0.5em}
{\rm (i)} Fix $i\le k$ and $\tau_{1},\tau_{2}\in \Tc$ such that $\sigma_{i}\le \tau_{1} \le \tau_{2}\le \sigma_{i+1}$ $\as$   Let us show that 
$$
\Xb_{\tau_{1}}\ge \Ec_{\tau_{1},\tau_{2}}[\Xb_{\tau_{2}}].
$$
Indeed, since $X\ge X_{+}$  and $I\ge 0$ by (b), \hSLD~implies that 
\begin{align*}
\Ec_{\tau_{1},\tau_{2}}[\Xb_{\tau_{2}}]&=\Ec_{\tau_{1},\tau_{2}}\left[X_{\tau_{2}}+I_{\sigma_{i}}+(X_{\sigma_{i}}-X_{\sigma_{i}+})\1_{\{\sigma_{i}<\tau_{2}\}}\right]\\
&\le \Ec_{\tau_{1},\tau_{2}}[X_{\tau_{2}}]+I_{\sigma_{i}}+(X_{\sigma_{i}}-X_{\sigma_{i}+})\1_{\{\sigma_{i}<\tau_{2}\}}.
\end{align*}
On the other hand,  it follows from (b) that 
$\Ec_{\tau_{1},\tau_{2}}[X_{\tau_{2}}]\le X_{\tau_{1}+}$. Hence, 
$$
	\Ec_{\tau_{1},\tau_{2}}[\Xb_{\tau_{2}}]\le  X_{\tau_{1}+}+I_{\sigma_{i}}+(X_{\sigma_{i}}-X_{\sigma_{i}+})\1_{\{\sigma_{i}<\tau_{2}\}}=X_{\tau_{1}}+I_{\tau_{1}}=\Xb_{\tau_{1}}.
$$

{\rm (ii)} In view of \hTC(b), the result of (i)  implies in particular that $\Xb_{\tau_{1}}\ge \Ec_{\tau_{1},\tau_{2}}[\Xb_{\tau_{2}}]$ for any  $(\tau_{1}, \tau_{2})\in \Tcd$ such that $\sigma_{i}\le \tau_{1}\le \sigma_{i+1}$ and $\sigma_{j}\le \tau_{2}\le \sigma_{j+1}$  a.s., for some $i\le j\le k$.  

\vspace{0.5em}
{\rm (iii)} Given $\tau \in \Tc$, we next show by induction that 
$$
\Xb_{\sigma_{i}}\ge \Ec_{\sigma_{i},  \tau}[\Xb_{  \tau}] \;\mbox{ on } \{\sigma_{i}\le \tau\}, \;\forall\;i\le k.
$$
For $i=k$, this follows from \hTC(c)~and (i).  
Assume that it is true for $1\le i+1\le k$. Then,  on $ \{\sigma_{i}\le \tau\}$, 
\begin{equation*}
\Xb_{\sigma_{i}}\ge \Ec_{\sigma_{i}, \tau\wedge \sigma_{i+1}}\left[\Xb_{\tau\wedge \sigma_{i+1}}\right]=  \Ec_{\sigma_{i}, \tau\wedge \sigma_{i+1}}\left[\Xb_{\tau}\1_{\tau \le \sigma_{i+1}}+\Xb_{\sigma_{i+1}}\1_{\tau > \sigma_{i+1}}\right].
\end{equation*}
But, by (a) and (c) of \hTC~and the induction hypothesis, we deduce immediately 
$$\Xb_{\tau}\1_{\tau \le \sigma_{i+1}}= \Ec_{\tau\wedge \sigma_{i+1},\tau}\left[\Xb_{\tau}\right]\1_{\tau \le \sigma_{i+1}}\ \text{and}\ 
\Xb_{\sigma_{i+1}}\1_{\tau > \sigma_{i+1}}\ge \Ec_{\tau\wedge \sigma_{i+1},\tau}\left[\Xb_{\tau}\right]\1_{\tau > \sigma_{i+1}}.$$
It remains to appeal to \hSLD~to deduce that $\Xb_{\sigma_{i}}\ge \Ec_{\sigma_{i}, \tau\wedge \sigma_{i+1}}\circ \Ec_{\tau\wedge \sigma_{i+1},\tau}[\Xb_{\tau}]$, on $ \{\sigma_{i}\le \tau\}$,  and to conclude by \hTC(b).

\vspace{0.5em}
{\rm (iv)} We are in position to conclude this step. Fix $(\tau_{1},\tau_{2})\in \Tcd$.   Set $\tilde \tau^{i}_{1}:=(\tau_{1}\vee \sigma_{i})\wedge \sigma_{i+1}$. Then,  \hTC(c) implies that 
$\Ec_{\tau_{1},\sigma_{i+1}\wedge \tau_{2}}[\Xb_{\tau_{2}}]= \Ec_{\tilde \tau^{i}_{1},\sigma_{i+1}\wedge \tau_{2}}[\Xb_{\sigma_{i+1}\wedge \tau_{2}}]$ on $\{\sigma_{i}\le \tau_{1}\le \sigma_{i+1}\}$. But it follows from (iii), and the same arguments as above, that 
\begin{align*}
\Ec_{\tilde \tau^{i}_{1},\sigma_{i+1}\wedge \tau_{2}}\left[\Xb_{\sigma_{i+1}\wedge \tau_{2}}\right]
&=
\Ec_{\tilde \tau^{i}_{1},\sigma_{i+1}\wedge \tau_{2}}\left[\Xb_{\tau_{2}}\1_{\tau_{2} \le \sigma_{i+1}}+\Xb_{\sigma_{i+1}}\1_{\tau_{2} > \sigma_{i+1}}\right]\\
&\ge \Ec_{\tilde \tau^{i}_{1},\sigma_{i+1}\wedge \tau_{2}}\left[\Ec_{\tau_{2}\wedge \sigma_{i+1},\tau_{2}}\left[\Xb_{\tau_{2}}\right]\right]\\
&= \Ec_{\tilde \tau^{i}_{1},\tau_{2}}\left[\Xb_{\tau_{2}}\right].
\end{align*}
Recalling the result of (i), we conclude that,  on $\{\sigma_{i}\le \tau_{1}\le \sigma_{i+1}\}$, 
$$\Xb_{\tau_{1}}\ge \Ec_{\tau_{1},\sigma_{i+1}\wedge \tau_{2}}\left[\Xb_{\sigma_{i+1}\wedge \tau_{2}}\right]\ge \Ec_{\tau_{1},\tau_{2}}\left[\Xb_{\tau_{2}}\right].$$ Since $\cup_{i=0}^{k} \{\sigma_{i}\le \tau_{1}\le \sigma_{i+1}\} = \Omega$, this concludes the proof of this step. 

\vspace{0.5em}
\noindent{\rm (d) \sl We now provide a bound on $I^{\eps,k}_T$ defined by \eqref{eq:IX_eps_k}.} 

\vspace{0.5em}
\noindent Let $(\sigma_{i})_{i\le k}$ be as in (c) associated to the parameter $(\eps, k)$. 
We first prove by induction that 
$$
\Ec_{\sigma_{i},T}[I^{\eps,k}_{T}]
~\le~ 
I_{\sigma_{i}}+X_{\sigma_{i}}+\E^{\Q_{i}}_{\sigma_{i}}[X_{T}^{-}],\ i\le k+1,
$$
in which  $\Q_{i}\in \Qc$. The result is true for $i=k+1$, recall our convention $\sigma_{k+1}=T$ and \hTC(a). 
Let us assume that it holds for some $i+1\le k+1$. Then, by \hTC(a)-(b)~and \hSLD~combined with (b), 
\begin{align*}
\Ec_{\sigma_{i},T}\left[I^{\eps,k}_{T}\right]&=\Ec_{\sigma_{i},\sigma_{i+1}}\circ \Ec_{\sigma_{i+1},T} \left[I^{\eps,k}_{T}\right]
\\
&\le \Ec_{\sigma_{i},\sigma_{i+1}}\left[I^{\eps,k}_{\sigma_{i+1}}+X_{\sigma_{i+1}}+\E^{\Q_{i+1}}_{\sigma_{i+1}}\left[X_{T}^{-}\right]\right] \\
&=\Ec_{\sigma_{i},\sigma_{i+1}}\left[I^{\eps,k}_{\sigma_{i}}+X_{\sigma_{i}}-X_{\sigma_{i}+}+X_{\sigma_{i+1}}+\E^{ \Q_{i+1}}_{\sigma_{i+1}}[X_{T}^{-}]\right] 
\\
&\le I^{\eps,k}_{\sigma_{i}}+X_{\sigma_{i}}-X_{\sigma_{i}+}+\Ec_{\sigma_{i},\sigma_{i+1}}\left[X_{\sigma_{i+1}}\right]+\E^{\tilde  \Q_{i}}_{\sigma_{i}}\left[\E^{\Q_{i+1}}_{\sigma_{i+1}}[X_{T}^{-}]\right]\\
&\le  I^{\eps,k}_{\sigma_{i}}+X_{\sigma_{i}}+\E^{\tilde \Q_{i}}_{\sigma_{i}}\left[\E^{\Q_{i+1}}_{\sigma_{i+1}}[X_{T}^{-}]\right] ,
\end{align*}
in which   $\tilde \Q_{i}\in \Qc$. Then, our induction claim follows for $i$ by composing $\tilde \Q_{i}$ and $\Q_{i+1}$ in an obvious way.
Recalling our convention $\sigma_{0}=0$, this implies that 
$
\Ec_{0,T}[I^{\eps,k}_{T}]\le X_{0}+\E^{\Q_{0}}[X_{T}^{-}]
$, from which \hSLD~provides the estimate 
$$
L^{-1}\E^{\Q}\left[I^{\eps,k}_{T}\right]
~\le~
\Ec_{0,T}\left[I^{\eps,k}_{T}\right ]-\Ec_{0,T}[0]
~\le~ 
X_{0}+\E^{\Q_{0}}[X_{T}^{-}]-\Ec_{0,T}[0],
$$
in which $\Q\sim \P$ is such that  $\E^{\Q}[|d\P/d\Q|^{q}]\le L$.
Since $p$ and $q$ are conjugate, it remains to use H\"older's inequality to deduce that 
\begin{align}\label{eq: borne I}
\E\left[ \big(I^{\eps,k}_{T} \big)^{\frac1p}\right]^{p}\le C_{L}\left(1+|X_{0}|+\E[(X_{T}^{-})^{p}]^{\frac1p}+|\Ec_{0,T}[0]|\right),
\end{align}
for some $C_{L}>0$ which only depends on $L$.

\vspace{0.5em}
\noindent{(e) \sl  We now extend the bound \reff{eq: borne I} to the general case.} 

\vspace{0.5em}
\noindent
Notice that the r.h.s. of \eqref{eq: borne I} does not depend on $\eps$ nor $k$, so we can first send $k$ to $\infty$ and then $\eps$ to $0$ and apply the monotone convergence theorem,
to obtain that
$$
\E\left[ \big(I_{T} \big)^{\frac1p}\right]^{p}\le C_{L}\left(1+|X_{0}|+\E[(X_{T}^{-})^{p}]^{\frac1p}+|\Ec_{0,T}[0]|\right).
$$

\vspace{0.5em}
\noindent{(f) \sl  It remains to show that $\Xb :=X+I $ is a local $\Ec$-supermartingale.} 

\vspace{0.5em}
\noindent
Recall that $I$ is defined in \reff{eq: def I}, and $(I^{\eps, k}, \Xb^{\eps, k})$ are defined in \eqref{eq:IX_eps_k}.
Let $\vartheta_{n}$ be the first time when $I\ge n$. 
Note that $(\vartheta_{n})_{n\ge 1}$ is a.s.~ increasing and converges to $\infty$, this follows from {\rm (e)}. We know from {\rm (c)} that $\Xb^{\eps,k}$ is a  $\Ec$-supermartingale. Hence, for $(\sigma,\tau)\in \Tcd$, we have 
$$
	\Xb^{\eps,k}_{\sigma\wedge \vartheta_{n}}\ge \Ec_{\sigma\wedge \vartheta_{n},\tau\wedge \vartheta_{n}}\left[\Xb^{\eps,k}_{\tau\wedge \vartheta_{n}}\right].
$$
But $\Xb^{\eps,k}_{\vartheta}$ $\uparrow$ $\Xb_{\vartheta}$ $\as$~for any stopping time $\vartheta$, when one let $k$ first go to $\infty$ and then $\eps$ to $0$.  Since $\Xb_{\tau\wedge \vartheta_{n}} \in \Lb^{p}(\Fc_{\tau})$, by definition of $(\vartheta_{n})_{n\ge 1}$ and the fact that $X\in {\bf X}^{p}_r$, \hS(c)~implies that 
$$
	\Xb_{\sigma\wedge \vartheta_{n}}\ge \Ec_{\sigma\wedge \vartheta_{n},\tau\wedge \vartheta_{n}}\left[\Xb_{\tau\wedge \vartheta_{n}}\right],
$$
which concludes the proof. \ep

\section{Doob-Meyer-Mertens decomposition of $g$-supermar-tingale systems}
\label{sec: doob meyer bsde}

We now specialize to the context of $g$-expectations introduced by Peng \cite{peng1997} 
(notice however that we consider a slightly more general version).
The object is to provide a Doob-Meyer-Mertens decomposition of $g$-supermartingale systems without c\`adl\`ag conditions. This is our Theorem \ref{thm: doob-mertens decompo bsde} below. 

 \vspace{0.5em}
\noindent We assume that the space $(\Omega,\Fc,\P)$ carries a $d$-dimensional Brownian motion $W$, adapted to the filtration $\F$, which may be strictly larger than the natural (completed) filtration of $W$. 
Recall that $\F$ satisfies the usual conditions. 

\subsection{$g$-expectation and Doob-Meyer decomposition}\label{sec:mertens}

Fix some $p>1$. Let $g:(\omega,t,y,z)\in \Omega\x [0,T]\x \R\x \R^{d}\longmapsto g_{t}(\omega,y,z)\in \R$ be such that $(g_{t}(\cdot,y,z))_{t\le T}$ is $\F$-progressively measurable for every $(y,z)\in\R\times\R^d$ and 
\begin{eqnarray} \label{eq:LipCond_g}
&|g_{t}(\omega,y,z)-g_{t}(\omega,y',z')|\le L_{g}(|y-y'|+|z-z'|),&  \label{eq: g Lip} \\
&\;\forall\;(y,z),(y',z')\in \R\x \R^{d}, \;\mbox{ for } dt\x d\P-{\rm a.e.}~(t,\omega)\in [0,T]\x \Omega, &\nonumber
\end{eqnarray}
for some constant number $L_{g}>0$. We also assume that $(g_{t}(\omega,0,0))_{t\le T}$ satisfies the following integrability condition
\begin{equation}\label{eq:integ}
\E\left[\int_0^T|g_t(0,0)|^pdt\right]<\infty.
\end{equation}
In the following, we most of the time omit the argument $\omega$ in $g$.

\vspace{0.5em}
\noindent Given $(\sigma,\tau)\in \Tcd$ and $\xi \in \Lb^{p}(\Fc_{\tau})$, we set $\mathcal E^{g}_{\sigma,\tau}[\xi]:=Y_{\sigma}$ in which $(Y,Z,N)$ is the unique solution of 
\begin{equation}\label{eq: BSDEweak}
Y_{t}=\xi+\int_{t\wedge \tau}^{\tau} g_{s}(Y_{s},Z_{s})ds -\int_{t\wedge \tau}^{\tau} Z_{s}  \cdot dW_{s}-\int_{t\wedge\tau}^\tau dN_s, \;t\le T,
\end{equation}
such that $(Y,Z) \in \Sb^{p} \x \Hb^{p}$ and $N$ is a c\`adl\`ag $\mathbb F$-martingale orthogonal to $W$ in the sense that the bracket $[ W,N]$ is null, $\mathbb P$-a.s., and such that
$$\E\left[[N]_T^{\frac p2}\right]<+\infty.$$

\noindent The wellposedness of this equation follows from  {\cite[Thm 4.1]{BPTZ}, see also \cite[Thm 2]{kp} and  \cite{klimsiak} or  \cite[Prop. A.1]{pt} for the case $p=2$.} We also remind the reader that the introduction of the orthogonal martingale $N$ in the definition of the solution is necessary, since the martingale predictable representation property may not hold with a general filtration $\mathbb F$. The map ${\cal E}^g$ is usually called the $g$-expectation operator.

\vspace{0.5em}
\noindent 
We define $\Ec^{g}$-supermartingales, also called $g$-supermatingales, as in the previous section, for $\Ec=\Ec^{g}$, 
i.e. $X$ is a $\Ec^g$-supermatingale iff $X\in {\bf X}^{p}$ and  $X_{\sigma}\ge \Ec^g_{\sigma,\tau}[X_{\tau}]$ $\as$~for all $(\sigma,\tau)\in \Tcd$. 
For c\`adl\`ag $\Ec^{g}$-supermartingales, we have the following classical Doob-Meyer decomposition, {which is a consequence of the well-posedness of a corresponding reflected backward stochastic differential equation. Its proof is provided in the Appendix}  
(see also Peng \cite[Thm. 3.3]{peng1999monotonic} in the case of a Brownian filtration).\footnote{{ We emphasize that we work with a filtration which is not assumed to be quasi-left continuous, a case which, as far as we know, has never been covered in the literature. The main technical results needed to establish Proposition \ref{prop:doob-general} are given in our companion paper \cite{BPTZ}.}
} 
\vspace{0.5em}
\begin{Proposition}\label{prop:doob-general}
	Let $X\in {\bf X}^p$ be a c\`adl\`ag $\Ec^{g}$-supermartingale. Then there exists $Z\in \Hb^{p}$,  a c\`adl\`ag process $A \in \Ab^{p}$
	and a c\`adl\`ag martingale $N$, orthogonal to $W$, satisfying $\E[[N]_T^{p/2}]<\infty$, such that 
	\b*
		X_{\sigma}=X_{\tau}+\int_{\sigma}^{\tau} g_{s}(X_{s},Z_{s})ds +A_{\tau}-A_{\sigma}-\int_{\sigma}^{\tau} Z_{s} \cdot dW_{s}-\int^\tau_\sigma dN_s,
	\e*
	for all $(\sigma,\tau)\in \Tcd$. Moreover, this decomposition is unique.
\end{Proposition}

\proof See Appendix \ref{subsec:RBSDE_DoobMeyer}.\ep

\vspace{2mm}

\subsection{Time consistence and  regularity  of $g$-expectations}

We now verify that the conditions of Lemma \ref{lem: X+I supermart abstract} apply to ${\cal E}^g$.

\begin{Proposition}\label{eq: prop Hyp ok if g decreasing2} Assume that $y\longmapsto g_{t}(\omega,y,z)$ is non-increasing for all $z\in \R$, for $dt\x d\P-{\rm a.e.}~(t,\omega)\in [0,T]\x \Omega$. Then, Assumptions \hTC, \hS~and \hSLD~hold for $\Ec^{g}$.
\end{Proposition}

\proof
First, notice that since $W$ is actually continuous, we not only have $[W,N]=0,\ {\rm a.s.}$, but also
$$\langle W, N\rangle = \langle W, N^{c}\rangle=\langle W, N^{d}\rangle=0,\ {\rm a.s.},$$
where $N^{c}$ (resp. $N^{d}$) is the continuous (resp. purely discontinuous) martingale part of $N$.
Then \hTC~follows from the definition of $\mathcal E^{g}$ and the uniqueness of a solution. 
The stability properties \hS(b) and (c)  follow from the path continuity of the $Y$ component of the solution of  \eqref{eq: BSDEweak} and the standard estimates given in {\cite[Thm~2.1 and Thm~4.1]{BPTZ}, see also  \cite[Prop.~3]{kp} for the case where the filtration is quasi left-continuous}\footnote{Notice that these estimates can be readily extended to the difference of two solutions of BSDEs, since, as pointed out in the proof of  \cite[Thm~4.2]{briand2003p}, such a difference is itself the solution to a BSDE. }.

\vspace{0.5em}
\noindent The fact that \hSLD~holds is a consequence of the usual linearization argument. Let $(Y,Z,N)$ and $(Y',Z',N')$ be the solutions of \reff{eq: BSDEweak} with terminal conditions $\xi$ and $\xi'$. Then, since $g$ is uniformly Lipschitz continuous, there exist two processes $\lambda$ and $\eta$, which are   $\F$-progressively measurable,  such that
$$ g_s(Y_s,Z_s)-g_s(Y'_s,Z'_s)=\lambda_s \left(Y_s-Y'_s\right)+\eta_s\cdot \left(Z_s-Z'_s\right),\ ds\times d\P-{\rm a.e.}$$
These two processes are  bounded by $L_{g}$ for $dt\x d\P-{\rm a.e.}$ $(t,\omega)\in [0,T]\x \Omega$, as a consequence of \reff{eq: g Lip}. Moreover, $\lambda\le 0$ since $g$ is non-increasing in $y$.

\vspace{0.5em}
\noindent Then, for any $0\leq t\leq s\leq T$, let us define the following continuous, positive and $\F$-progressively measurable process
$$A_{t,s}:=\exp\left(\int_t^s\lambda_udu-\int_t^s\eta_u \cdot dW_u-\frac12\int_t^s|\eta_u|^2d{u}\right).$$
By applying It\^o's formula, we deduce classically (see \cite[Lem.~9]{kp}) that
$$Y_\sigma-Y'_\sigma=\mathbb E_\sigma\left[A_{\sigma,\tau}(\xi-\xi')\right],$$
which is nothing else but Assumption $\hSLD$ by Girsanov's theorem (recall that $\lambda\le 0$ and that $\lambda$ and $\eta$ are bounded  by $L_{g}$, i.e.~it suffices to consider $\Qc$ as the collection of measures with density with respect to $\P$ given by an exponential of Doléans-Dade of the above form with $\eta$ bounded by $L_{g}$). 

\vspace{0.5em}
\noindent Finally, the condition \hS(a) follows from a similar linearization argument.
Let $s \in [0,T)$ and $\xi \in \Lb^0(\Fc_s)$, $s_n \searrow s$ and $(\xi_n)_{n \ge 1}$ be such that $\xi_n \in \Lb^p(\Fc_{s_n})$ for each $n$, $(\xi^-_n)_{n \ge 1}$ is bounded in $\Lb^p$ and $\xi_n \to \xi$ as $n \to \infty$. One has
$$\Ecg_{s,s_{n}}[\xi_{n}]\ge  \E_{s}\left[A_{n}\left(\xi_{n}-C\int_{s}^{s_{n}}|g_{s}(0,0)|ds\right)\right],$$ for a sequence $(A_{n})_{n \ge 1}$ bounded in any $\Lb^{p'}$, $p'\ge 1$, which converges {\rm a.s.}~to $1$, and some $C\ge 1$ independent on $n$. Since $\big( \xi_{n}^{-},$ $\int_{s}^{s_{n}}|g_{s}(0,0)|ds \big)_{n \ge 1}$ is bounded in $\Lb^{p}$, and $p>1$, the negative part of term in the above expectation  is uniformly integrable, and we can apply Fatou's Lemma to conclude the proof.
\ep

\begin{Remark}\label{rem: reinforcement S}{\rm One easily checks that $X_{\sigma+}\ge \Ec^{g}_{\sigma,\tau}[X_{\tau+}]$ for $(\sigma,\tau)\in \Tcd$, whenever $X$ is a $\Ec^{g}$-supermartingale. Again, this follows from the path continuity of the $Y$ component of the solution of  \reff{eq: BSDEweak} and the estimates of {\cite[Rem 4.1]{BPTZ}}.}\end{Remark}

\begin{Corollary}\label{cor: application BSDE g monotone} Assume that $y\longmapsto g_{t}(\omega,y,z)$ is non-increasing for all $z\in \R^{d}$, for $dt\x d\P-{\rm a.e.}~(t,\omega)\in [0,T]\x \Omega$.   Let $X\in {\bf X}^{p}_{r}$ be an $\Ec^{g}$-supermartingale. Define the process $I$ by 
\begin{equation}\label{eq: def I2}
I_{t}:=\sum_{s<t}  (X_{s}-X_{s+}),\;t\le T.
\end{equation}
Then, $I$ is a non-decreasing and left-continuous process satisfying $I_{T}\in \Lb^{\frac1p}$. Moreover, $\Xb :=X+I$ is a right-continuous local $\Ec^{g}$-supermatingale.
\end{Corollary}

\proof This is an immediate consequence of Lemma \ref{lem: X+I supermart abstract} and Proposition \ref{eq: prop Hyp ok if g decreasing2} if $(X_{t}^{-})_{t\le T}$ is bounded in $\Lb^{p}$. But this follows from the fact that $X^{-}\le \Ecg_{\cdot,T}[X_{T}]^{-}\in \Sb^{p}$.
\ep

\vspace{0.5em}
\noindent For later use, let us a provide another version in which the monotonicity of $g$ in $y$ is not used anymore. The price to pay is that the $I$ process defined below may not be non-decreasing anymore.

\begin{Corollary}\label{cor: application BSDE g non monotone}  Let $X\in {\bf X}^{p}_{r}$ be an $\Ec^{g}$-supermartingale. Then, $X_t \ge X_{t+}$ for all $t \in [0,T)$. Define the process $I$ by 
\begin{equation}\label{eq: def I3}
I_{t}:=\sum_{s<t} e^{L_g({s}-{t})}(X_{s}-X_{s+}),\;t\le T.
\end{equation}
Then,  $I_{T}\in \Lb^{\frac1p}$, $I$  is left-continuous. Moreover, $\Xb:=X+I$ is a right-continuous local $\Ec^{g}$-supermatingale.
\end{Corollary}

\proof It follows from Corollary \ref{cor: application BSDE g monotone} that the result holds if $g$ is non-increasing in its $y$-variable. 
On the other hand, it is immediate to check that $\zeta$ is an $\Ec^{g}$-supermartingale if and only if $\tilde \zeta$ is an $\Ec^{\tilde g}$-supermartingale, with 
\b*
\tilde \zeta:=e^{L_gT} \zeta \;\mbox{ and }\;\tilde g_t(y,z):=e^{L_gt}g(ye^{-L_gt},ze^{-L_gt})-L_gy.
\e*
The map $\tilde g$ is now non-increasing in its $y$-component as a consequence of \reff{eq: g Lip}. Moreover, $\tilde g$ still satisfies \reff{eq: g Lip}, with the same constant $L_g$, and, by \reff{eq:integ}, $\tilde g(0,0)$ satisfies the  integrability condition needed to define the corresponding BSDE. Hence $\tilde X+\tilde I$ is a right-continuous $\Ec^{\tilde g}$-supermartingale, $\tilde I$ is non-decreasing and $\tilde I_{T} \in \Lb^{\frac{1}{p}}$, where we have defined 
$$
\tilde I_{t}:=\sum_{s<t} (\tilde X_{s}-\tilde X_{s+})=\sum_{s<t} e^{L_g{s}}( X_{s}- X_{s+}).
$$
Hence, $X+I=e^{-L_g\cdot}(\tilde X+\tilde I)$ is a $\Ec^{g}$-supermartingale, 
and $I_{T}\in \Lb^{\frac1p}$ since $\tilde I_{T} \in \Lb^{\frac{1}{p}}$.
\ep

\subsection{The Doob-Meyer-Mertens's decomposition for $\Ec^g$-supermar-tingales}

We are now in position to state the main result of this paper.
\vspace{2mm}

\noindent Let $S=\{S(\tau),\;\tau\in \Tc\}$ be a $\Tc$-system in the sense that for all $\tau,\tau'\in \Tc$
\begin{enumerate}[{\rm (i)}]
\item $S(\tau)\in \Lb^{0}(\Fc_{\tau})$,
\item  $S(\tau)=S(\tau')$ $\as$~on $\{\tau=\tau'\}$.
\end{enumerate}

\noindent If  $S(\tau)\in \Lb^{p}(\Fc_{\tau})$ for every $\tau \in \Tc$ and $S(\sigma)\ge \Ec^{g}_{\sigma,\tau}[S(\tau)]$ for all $(\sigma,\tau)\in \Tcd$, then we say that it is a $\Ec^{g}$-supermartingale system.

\begin{Theorem}[Mertens's decomposition]\label{thm: doob-mertens decompo bsde}  Let $S$ be a $\Ec^{g}$-supermatingale system s.t.~$\{S(\tau),$ $\tau \in \Tc\}$ is uniformly integrable, then there exists $X\in {\bf X}_{\ell r}^{p}$ such that $S(\sigma)=X_{\sigma}$ for all $\sigma \in \Tc$. 
If in addition,  $\esssup\{S(\tau)$ $\tau \in \Tc\}\in \Lb^{p}$, then  there exists $(Z,A)\in   \Hb^{p}\x \Ab^{p}$ and a c\`adl\`ag martingale $N$, orthogonal to $W$, satisfying $\E[[N]_T^{p/2}]<\infty$, such that 
\b*
S(\sigma)=X_{\sigma}=X_{\tau}+\int_{\sigma}^{\tau} g_{s}(X_{s},Z_{s})ds +A_{\tau}-A_{\sigma}-\int_{\sigma}^{\tau} Z_{s}  \cdot dW_{s}-\int_\sigma^\tau dN_s,
\e*
for all $(\sigma,\tau)\in \Tcd$. This decomposition is unique. 
\end{Theorem}

\proof {\rm (a)} Let us first prove that there exists an optional process $X\in {\bf X}^{p}$ such that $S(\sigma)=X_{\sigma}$ a.s.~for all $\sigma\in \Tc$. 
 {Since $S$ is  uniformly integrable,     \cite[Thm.~6 and Rem.~7 c)]{dellacherie1982problemes} imply that it suffices to show that }
$$
\E[S(\sigma)]\ge \liminf_{n\to \infty}\E[ S(\sigma_{n})],
$$
for all non-increasing sequence $(\sigma_{n})_{n\ge 1} \in \Tc_{\sigma}$ such that $\sigma_{n}\longrightarrow \sigma\in \Tc$, $\as$ By using a similar linearization argument as the one used in the proof of Proposition \ref{eq: prop Hyp ok if g decreasing2}, we can find $\F$-progressively measurable processes $\lambda^{n}$ and $\eta^{n}$ that are bounded by $L_{g}$ $dt\x d\P$-a.e.~and such that 
$$
S(\sigma)\ge \E_{\sigma}\left[H^{n}_{\sigma_{n}} \left(e^{\int_\sigma^{\sigma_n}\lambda^n_sds}S(\sigma_{n})+\int_{\sigma}^{\sigma_{n}}e^{\int_{\sigma}^{s} \lambda^{n}_{s}d{s}}g_{s}(0,0)d{s}\right)\right]
$$
where 
$$
H^{n}:=\exp\left(-\frac12\int_{\sigma}^{\cdot \vee \sigma } |\eta^{n}_{s}|^{2}d{s}+ \int_{\sigma}^{\cdot\vee \sigma} \eta^{n}_{s}  \cdot dW_{s} \right). 
$$
Then, 
\begin{align}
\E[S(\sigma)]\ge&\ \E[S(\sigma_{n})]+ \E\left[e^{\int_\sigma^{\sigma_n}\lambda^n_sd{s}}(H^{n}_{\sigma_{n}}-1) S(\sigma_{n})\right]\nonumber\\
&+ \E\left[H^{n}_{\sigma_{n}}\int_{\sigma}^{\sigma_{n}}e^{\int_{\sigma}^{s} \lambda^{n}_{s} d{s}}g_{s}(0,0)d{s} \right]\label{eq: borne inf  pour aggreg}\\
&+\E\left[\left(e^{\int_\sigma^{\sigma_n}\lambda^n_sd{s}}-1\right)S(\sigma_n)\right]. \nonumber
\end{align}
Note that
$$(H^n_{\sigma_n}-1)S(\sigma_n)\geq -(S(\sigma_n))^+-H^n_{\sigma_n}(S(\sigma_n))^-.$$
Since $S$ is uniformly integrable, so is $S^+$. Besides, we have by definition
$$S(\sigma_n)\geq \Ec^{g}_{\sigma_n,T}[S(T)].$$ 
But, once more it is clear that $\Ec^{g}_{\sigma_n,T}[S(T)]$ is bounded in $\Lb^p$, uniformly in $n$, {see \cite[Thm 4.1]{BPTZ}}. Since $H^n_{\sigma_n}$ has bounded (uniformly in $n$) moments of any order, 
de la Vall\'ee-Poussin criterion ensures that $H^nS^-$ is also uniformly integrable. Therefore, $\{[(H^n_{\sigma_n}-1)S(\sigma_n)]^{-},n\ge 1\}$ is uniformly integrable. Using the fact that $(\lambda^n,\eta^{n})_{n}$} is uniformly bounded by $L_g$, as well as \reff{eq:integ}, we can use Fatou's lemma in \reff{eq: borne inf  pour aggreg} to obtain that the second and the third terms on the right-hand side converges to $0$ as $n\longrightarrow \infty$.

\vspace{0.5em}
(b) The fact that $X$ has right- and left-limits, up to an evanescent set, follows from Lemma  \ref{lem: version limited a droite} stated below, since $X$ is an $\Ec^{g}$-supermartingale.

\vspace{0.5em}
(c) Let $I$ be defined as in Corollary \ref{cor: application BSDE g non monotone}  for $X$. Since $X+I$ is right-continuous, we can apply the Doob-Meyer decomposition of Proposition \ref{prop:doob-general} to $\Xb^{n}:=(X+I)_{\cdot\wedge \vartheta_{n}}$ where $\vartheta_{n}$ is the first time when $I\ge n$. There  exists  $(Z^{n},\bar A^{n})\in \Hb^{p}\x \Ab^{p}$ and a c\`adl\`ag martingale $N^n$, orthogonal to $W$, such that, for $(\sigma,\tau)\in \Tcd$, 
\begin{align*}
\Xb^{n}_{\sigma}&=\Xb^{n}_{\tau}+\int_{\vartheta_{n}\wedge \sigma}^{\vartheta_{n}\wedge \tau} g_{s}(\Xb^{n}_{s},Z^{n}_{s})ds+\bar A^{n}_{\tau}-\bar A^{n}_{\sigma} -\int_{\vartheta_{n}\wedge \sigma}^{\vartheta_{n}\wedge \tau} Z^{n}_{s} \cdot dW_{s}-\int_{\vartheta_{n}\wedge \sigma}^{\vartheta_{n}\wedge \tau}dN^n_s
\\
&= \Xb^{n}_{\tau}+\int_{\vartheta_{n}\wedge \sigma}^{\vartheta_{n}\wedge \tau} \{g_{s}(X_{s},Z^{n}_{s})+\eta_{s} I_{s}\} ds+\bar A^{n}_{\tau}-\bar A^{n}_{\sigma} -\int_{\vartheta_{n}\wedge \sigma}^{\vartheta_{n}\wedge \tau} Z^{n}_{s} \cdot dW_{s}\\
&\hspace{0.9em}-\int_{\vartheta_{n}\wedge \sigma}^{\vartheta_{n}\wedge \tau}dN^n_s,
\end{align*}
in which $\eta$ is a progressively measurable process bounded by $L_{g}$, $dt\x d\P$-a.e., as a consequence of \reff{eq: g Lip}. Set 
\begin{equation}\label{eq: def An}
A^{n}:=I_{\cdot \wedge \vartheta_{n}}+\bar A^{n}+\int_{0}^{\cdot \wedge \vartheta_{n}} \eta_{s} I_{s} ds,
\end{equation}
and observe that $(A^{n},Z^{n},N^n)=(A^{k},Z^{k},N^k)$ on $[\![0,\vartheta_{k}]\!]$ for $n\le k$, by uniqueness of the decomposition in Proposition \ref{prop:doob-general}. 
We can then define 
\begin{equation}\label{eq: def A par An}
(A,Z,N):=\1_{[\![0,\vartheta_{1}]\!]} (A^{1},Z^{1}, N^1)+ \sum_{n\ge 1} \1_{]\!]\vartheta_{n},\vartheta_{n+1}]\!]} (A^{n+1},Z^{n+1},N^{n+1}),
\end{equation}
so that 
\begin{align}\label{eq: dyna X avec A avant savoir est croissant}
X_{\sigma}&=X_{\tau}+\int_{ \sigma}^{  \tau} g_{s}(X_{s},Z_{s})ds+A_{\tau}- A_{\sigma} -\int_{ \sigma}^{ \tau} Z_{s} \cdot dW_{s}-\int_\sigma^\tau dN_s.
\end{align}
We claim  that $A$ is non-decreasing and that the above decomposition is unique. 
The fact that  $(Z, A,[N]_T)\in \Hb^{p}\x \Ab^{p}\x\Lb^{p/2}$  {then follows from    \cite[Proposition 2.1]{BPTZ}}.

\vspace{0.5em}
\noindent Let us now prove our claim. Define $\tilde X$ and $\tilde g$  by 
\b*
\tilde X:=e^{L_g\cdot }X \;\mbox{ and }\;\tilde g_t(y,z):=e^{L_gt}g(ye^{-L_gt},ze^{-t})-L_gy.
\e*
Then, $\tilde X$ is a $\Ec^{\tilde g}$-supermartingale, and so is its right-limits process $\tilde X^{+}:=\tilde X_{+}$, as a consequence of Remark \ref{rem: reinforcement S}, recall Lemma \ref{lem: version limited a droite} below. Applying Proposition \ref{prop:doob-general}, we can find a right-continuous non-decreasing process $\tilde A\in\Ab^p$, $\tilde Z\in \Hb^{p}$ and a c\`adl\`ag martingale $\tilde N$, orthogonal to $W$, such that 
\begin{align*}
\tilde X^{+}_{\sigma}&=\tilde X^{+}_{\tau}+\int_{ \sigma}^{  \tau} \tilde g_{s}(\tilde X^{+}_{s},\tilde Z_{s})ds+\tilde A_{\tau}-\tilde A_{\sigma} -\int_{ \sigma}^{ \tau} \tilde Z_{s} \cdot dW_{s}-\int_{ \sigma}^{ \tau}d\tilde N_s
\end{align*}
for all $(\sigma,\tau) \in \Tcd$. This decomposition is unique. On the other hand, \reff{eq: dyna X avec A avant savoir est croissant} implies that 
\begin{align*}
\tilde X^{+}_{\sigma}&=\tilde X^{+}_{\tau}+\int_{ \sigma}^{  \tau} \tilde g_{s}(\tilde X^{+}_{s},e^{L_{g}s} Z_{s})ds +\tilde B_{\tau}-\tilde B_{\sigma}-\int_{ \sigma}^{ \tau} e^{L_{g}s} Z_{s} \cdot dW_{s}-\int_\sigma^\tau e^{L_{g}{s}} dN_s,
\end{align*}
in which 
\begin{align*}
\tilde B_{\tau}-\tilde B_{\sigma}:=&\ \int_{ \sigma}^{  \tau} (\tilde g_{s}(\tilde X_{s},e^{L_{g}s} Z_{s})-\tilde g_{s}(\tilde X^{+}_{s},e^{L_{g}{s}} Z_{s}))ds\\
&+\int_{\sigma}^{\tau} e^{L_{g}{s}}dA_{s}+e^{L_{g}{\tau}}(A_{\tau+}-A_{\tau})-e^{L_{g}{\sigma}}(A_{\sigma+}-A_{\sigma}).
\end{align*}
Hence, $\tilde B=\tilde A$ is non-decreasing.  But, since 
$(\tilde g(\tilde X,e^{L_{g}\cdot} Z)-\tilde g(\tilde X^{+},e^{L_{g}\cdot} Z))\le 0$
as a consequence of Corollary \ref{cor: application BSDE g non monotone} (namely $\tilde X\ge \tilde X^{+}$) and the fact that $\tilde g$ is non-increasing in its first component,  we must have that the continuous part of $\int_{0}^{\cdot} e^{L_{g}{s}}dA_{s}$ is non-decreasing, and so must be the continuous part of $A$.  We now deduce from the definition of $I$ in \reff{eq: def I3} and \reff{eq: def An}-\reff{eq: def A par An} that $A$ can only decrease in a continuous manner, recall that $\bar A_{n}$ is non-decreasing. Hence, $A$ is non-decreasing. 
The fact that the decomposition is unique comes from the uniqueness of the decomposition for $\tilde X^{+}$.
\ep

\subsection{Remarks}

The framework of this section corresponds to the case where the BSDEs are driven by a continuous martingale $M$, whose quadratic variation is absolutely continuous with respect to the Lebesgue measure, and with an invertible density. Extensions to the context of \cite{brianddelyon2}, see also    \cite{elkarh}, \cite{klimsiak} or \cite{cohen}, would be of interest. Similarly, one could certainly consider BSDEs with jumps, generators with quadratic growth, obstacles, stochastic Lipschitz conditions, etc.
We have chosen to work in a simpler setting so as not to drown our arguments with unneeded technicalities, and to focus on the novelty of our approach.
  
  \vspace{0.5em}
\noindent However, the case  $p=1$ can not be treated by the same technics, in particular we can not appeal to the classical linearization procedure. It would also require a reinforcement of the condition \reff{eq:LipCond_g}, see \cite{briand2003p}.

\def\u{{\rm u}}
\def\v{{\rm v}}
\def\dO{{\rm d}_{\Oc}}
\def\Bc{\mathcal{B}}
\def\Pc{\mathcal{P}}

\section{Applications}

We now consider two problems studied in the recent literature, which are solved with sophisticated arguments under technical conditions.
Using Theorem \ref{thm: doob-mertens decompo bsde}, we can solve these problems in a very general context with quite simple arguments.

\subsection{Optional decomposition of $g$-supermartingale systems}

\noindent We are still in the context of the previous section, with the slight modification   that,  instead of  the Brownian motion $W$, we consider a continuous $(\P,\F)$-martingale $M$ of the form: $M_t = \int_{0}^t \alpha^{\top}_{s}  dW_{s}$, in which $\alpha$ is a $\R^{d\times d}$-bounded predictable process with bounded inverse.
Recall that $\F$ satisfies the usual conditions.

\vspace{0.5em}
\noindent Let $S=\{S(\tau),\;\tau\in \Tc\}$ be a $\Tc$-system, $g$ be as in Section \ref{sec: doob meyer bsde} such that \reff{eq:LipCond_g} and \reff{eq:integ} hold. 
Let $ \Mc_0$ denote the set of probability measures $\mathbb Q$ on $(\Omega, \mathcal F)$ which are equivalent to $\P$ and such that $M$ is a $(\Q,\F)$-martingale. We then say that a $\mathcal T$-system $S$ is a $\mathcal E^{g}$-supermartingale system under some $\Q\in\Mc_0$ if $S(\tau)\in \Lb^{p}(\Q)$ for all $\tau\in\mathcal T$ and $S(\sigma)\ge \Ec^{\Q,g}_{\sigma,\tau}[S(\tau)]$ for all $(\sigma,\tau)\in \Tcd$, where, with $(\sigma,\tau)\in \Tcd$ and $\xi \in \Lb^{p}(\Fc_{\tau},\Q)$, we set $\mathcal E^{\Q,g}_{\sigma,\tau}[\xi]:=Y_{\sigma}$, with $(Y,Z,N)$ the unique solution of 
\begin{align*}
Y_{t}&=\xi+\int_{t\wedge \tau}^{\tau} g_{s}(Y_{s},Z_{s})ds -\int_{t\wedge \tau}^{\tau} Z_{s}  \cdot dM_{s}-\int_{t\wedge\tau}^\tau dN_s,\\
&=\xi+\int_{t\wedge \tau}^{\tau} g_{s}(Y_{s},\alpha^{-1}_{s}\alpha_{s}Z_{s})ds -\int_{t\wedge \tau}^{\tau} \alpha_{s}Z_{s}  \cdot dW_{s}-\int_{t\wedge\tau}^\tau dN_s,\; \;t\le T,
\end{align*}
such that $Y\in \Sb^{p}(\Q)$, $Z$ belongs to $\Hb^{p}(\Q)$ and $N$ is a c\`adl\`ag $(\mathbb F,\Q)$-martingale orthogonal to $M$, and such that
$$\E^\Q\left[[N]_T^{p/2}\right]<+\infty.$$
The spaces $\Sb^{p}(\Q)$ and $\Hb^{p}(\Q)$ are defined as $\Sb^p$ and $\Hb^p$, but with $\Q$ instead of $\P$.

\vspace{0.5em}
\noindent The main result of this section is the following optional type decomposition (see e.g. \cite{elkarq, Kramkov, fk98}).

\begin{Theorem}[Optional decomposition]\label{thm: optional decompo bsde} If for any $\Q\in\Mc_0$, $S$ is a $\Ec^{\Q,g}$-super-martingale system which is $\Q$-uniformly integrable s.t.  $\esssup\{|S(\tau)|,$ $\tau\in\mathcal T\}\in \Lb^{p}(\Q)$, then there exists $(X,Z)\in {\bf X}^{p}_{\ell r}\x \Hb^{p}$ such that $S(\sigma)=X_{\sigma}$ for all $\sigma\in \Tc$, and
\b*
X_{\cdot}+\int_{0}^{\cdot} g_{s}(X_{s},Z_{s})ds -\int_{0}^{\cdot} Z_{s} \cdot dM_{s}\text{ is non-increasing, ${\rm a.s.}$}
\e*
\end{Theorem}

\proof The existence of the process $X\in {\bf X}^{p}_{\ell r}$ such that $S(\sigma)=X_{\sigma}$ for all $\sigma\in \Tc$ follows from Theorem \ref{thm: doob-mertens decompo bsde}. Fix then some $\Q\in\Mc_0$. Using Theorem \ref{thm: doob-mertens decompo bsde}, we deduce the existence of $(Z^\Q,A^\Q)\in \Hb^{p}(\Q)\x \Ab^{p}(\Q)$ and of a $\Q$-martingale $N^\Q$ orthogonal to $M$ such that $ \P-{\rm a.s}. \ \text{(or $\Q-{\rm a.s}.$)}$
\begin{equation*}
X_{\sigma}=X_{\tau}+\int_{\sigma}^{\tau} g_{s}(X_{s},Z^\Q_{s})ds +A^\Q_{\tau}-A^\Q_{\sigma}-\int_{\sigma}^{\tau} Z^\Q_{s}\cdot  dM_{s}-\int_t^TdN_s^\Q,
\end{equation*}
for $(\sigma,\tau)\in \Tcd$.
Recall the definition of $I$ in Corollary \ref{cor: application BSDE g non monotone} and that $X+I$ is right-continuous. Then, 
\begin{equation}\label{eq: croche X+I}
[X+I,M]_\cdot=\int_0^\cdot   \alpha^{\top}_{s}\alpha_{s} Z^{\Q}_s ds,
\end{equation}
and  the family $(Z^\Q)_{\Q\in\Mc_0}$ can actually be aggregated into a universal predictable process $Z$, since $\alpha$ is invertible. Hence, we deduce that $X+\int_0^\cdot g_s(X_s,Z_s)ds$ is actually a supermartingale under any $\Q\in\Mc_0$, and we can apply the classical optional decomposition theorem (\cite[Thm.1]{fk98}) together with the classical Mertens's decomposition (\cite[T2 Lemme]{mertens1972theorie}) to deduce the existence of an $\F$-predictable process $\widetilde Z$ such that
	$$
		X_{\cdot}+\int_{0}^{\cdot} g_{s}(X_{s},Z_{s})ds -\int_{0}^{\cdot} \widetilde Z_{s} \cdot dM_{s}
		\text{ is non-increasing, $\P-{\rm a.s.}$}
	$$
Next, using \eqref{eq: croche X+I}, we obtain 
$Z=\widetilde Z$  $dt \times d\P$-{\rm a.e.},
which ends the proof.
\ep

\subsection{Dual formulation for minimal super-solutions of BSDEs with constraints on the gains process}

 In this section, we provide an application to the  dual representation for BSDEs with constraints.   We specialize to the situation where $\Omega$ is the canonical space of $\R^{d}$-valued continuous functions on $[0,T]$, starting at $0$, endowed with the Wiener measure $\P$. We let $\F^{\circ}=(\Fc^{\circ}_{t})_{t\le T}$ denote the raw filtration of the canonical process $\omega\longmapsto W(\omega)=\omega$, while $\F$  denotes its $\P$-augmentation. We also fix $p'>p>1$.
  
\vspace{0.5em}
\noindent We let $g$ be as in Section \ref{sec: doob meyer bsde} such that \reff{eq:LipCond_g} and \reff{eq:integ} hold for $p'$ and fix $\xi \in \Lb^{p'}$.
Further, let $\Oc = ( \Oc_{t}(\omega) )_{(t,\omega) \in [0,T] \x \Omega}$ be
a family of non-empty closed convex random subsets of $\R^{d}$,  
such that $\Oc$ is $\F^{\circ}$-progressively measurable in the sense of random sets (see e.g. Rockafellar~\cite{Rockafellar.76})
i.e.~$\{(s,\omega) \in [0,t] \x \Omega ~: \Oc_{s}(\omega) \cap O \neq \emptyset \} \in \Bc([0,t]) \otimes \Fc$  for all $t \in [0,T]$ and all closed $O \subseteq \R^d$.
In particular, it admits   a Castaing representation, see e.g.~\cite{Rockafellar.76}, which in turn ensures that  the support function defined by
 $$
 \delta_{t}(\omega,\cdot) : u\in \R^{d}\longmapsto \delta_{t}(\omega,u) := \sup\{u\cdot z, \; z\in \Oc_{t}(\omega)\}
 $$
is $\Fc^{\circ}_{t}\otimes \Bc([0,t])\otimes \Bc(\R^{d})/\Bc(\R^{d}\cup \{\infty\})$-measurable, for each $t \in[0, T]$.
 
\vspace{0.5em}
\noindent We consider solutions $(Y,Z,A)\in  {\bf X}^{p}_{\ell r} \x \Hb^{p}\x \Ab^{p}$ of 
\be\label{eq: bsde supersol}
Y=\xi+\int_{\cdot}^{ T} g_{s}(Y_{s},Z_{s})ds+A_{T}- A -\int_{ \cdot}^{ T} Z_{s}   \cdot dW_{s}, 
\ee
under the constraint 
 \be\label{eq: constraint z}
 Z\in \Oc, \;dt\x d\P-{\rm a.e.} 
 \ee
 
\vspace{0.5em}
\noindent We say that a solution $(Y,Z,A)\in  {\bf X}^{p}_{\ell r}\x \Hb^{p}\x \Ab^{p}$ of \reff{eq: bsde supersol}-\reff{eq: constraint z} is minimal if any other solution $(Y',Z',A')\in {\bf X}^{p}_{\ell r} \x \Hb^{p}\x \Ab^{p}$ is such that $Y_{\tau}\le Y'_{\tau}$ $\as$,~for any $\tau \in \Tc$.

 \vspace{0.5em}
\noindent The dual characterization relies on the following construction. 

\vspace{0.5em}
\noindent Let us also define $\Uc$ as the class of $\R^d$-valued, progressively measurable processes such that $ |\nu|+|\delta(\nu)|\le c$, $dt\x d\P$-a.e.,~for some $c\in \R$.
 Given $\nu \in \Uc$, we let $\P^{\nu}$ be the probability measure whose density with respect to $\P$ is given by the Dol\'eans-Dade exponential of $\int_{0}^{\cdot} \nu_{s} \cdot dW_{s}$, and denote by $W^{\nu}:=W-\int_{0}^{\cdot} \nu_{s} ds$ the corresponding $\P^{\nu}$-Brownian motion.     Then, given $\xi' \in \Lb^{p}(\Fc_{\tau},\P^{\nu})$, $\tau \in \Tc$, we define $\Ec^{\nu}_{\cdot,\tau}[\xi']$ as the $Y^{\nu}$-component of the solution $(Y^{\nu},Z^{\nu})\in \Sb^{p}(\P^{\nu})\x \Hb^{p}(\P^{\nu})$ of the BSDE
 $$
 Y^{\nu}=\xi'+\int_{\cdot}^{\tau} \left(g_{s}(Y^{\nu}_{s},Z^{\nu}_{s})-\delta_{s}(\nu_{s})\right)ds -\int_{\cdot}^{\tau} Z^{\nu}_{s} \cdot dW^{\nu}_{s}. 
 $$ 
In the above,  $\Sb^{p}(\P^{\nu})$ and $ \Hb^{p}(\P^{\nu})$ are defined as $\Sb^p$ and $\Hb^{p}$ but with respect to $\P^{\nu}$ in place of $\P$. 
 
 \begin{Theorem}\label{thm: dual bsde constraint} Define
 \be\label{eq: def S dual bsde constraint}
 S(\tau):=\esssup\left\{\Ec^{\nu}_{\tau,T}[\xi],\;\nu \in \Uc\right\},\;\;\tau \in \Tc.
 \ee
 Assume that  $\esssup\{|S(\tau)|, \;\tau\in\mathcal T\}\in \Lb^{p'}$ for some $p'>p$. Then, there exists   $X\in {\bf X}^{p}_{\ell r}$ such that $X_{\tau}=S(\tau)$ for all $\tau \in \Tc$, and $(Z,A)\in \Hb^{p}\x \Ab^{p}$ such that $(X,Z,A)$ is the minimal solution of \reff{eq: bsde supersol}-\reff{eq: constraint z}.
 \end{Theorem}
 
\noindent Before providing the proof of this result, let us comment it. This formulation is known since  \cite{cvitanic1998backward},  however it was proven only under strong assumptions. Although it should essentially be a consequence of the Doob-Meyer decomposition for $g$-supermatingales, the main difficulty comes from the fact that the family of controls in  $\Uc$ is not uniformly bounded. Hence, \reff{eq: def S dual bsde constraint} is a singular control problem for which the right-continuity of  $\tau\longmapsto S(\tau)$ is very difficult to establish, {\it a priori}, see  \cite{bouchard2014regularity} for a restrictive Markovian setting. This fact prevents us to apply the result of \cite{peng1999monotonic}. Theorem \ref{thm: doob-mertens decompo bsde} allows us to bypass this issue and provides a very simple proof.
 
 \vspace{0.5em}
 \noindent{\bf Proof of Theorem \ref{thm: dual bsde constraint}.} Let $(Y,Z,A)\in {\bf X}^{p}_{\ell r}\x \Hb^{p}\x \Ab^{p}$ be a solution of \reff{eq: bsde supersol}-\reff{eq: constraint z}. Then, for $(\sigma,\tau) \in \Tcd$, 
  \begin{align*}
 Y_{\sigma}=&\ Y_{\tau}+\int_{\sigma}^{\tau}  g_{s}(Y_{s},Z_{s}) ds +A_{\tau}-A_{\sigma} -\int_{\sigma}^{\tau} Z_{s}  \cdot  dW_{s}\\
 =&\ Y_{\tau}+\int_{\sigma}^{\tau} ( g_{s}(Y_{s},Z_{s}) - \nu_{s}\cdot Z_{s} )ds +A_{\tau}-A_{\sigma} -\int_{\sigma}^{\tau} Z_{s} \cdot  dW^{\nu}_{s}\\
 =&\ Y_{\tau}+\int_{\sigma}^{\tau} ( g_{s}(Y_{s},Z_{s}) - \delta_{s}(\nu_{s}) )ds +A_{\tau}-A_{\sigma} +\int_{\sigma}^{\tau} (\delta_{s}(\nu_{s}) -\nu_{s}\cdot Z_{s}) ds \\
 & -\int_{\sigma}^{\tau} Z_{s}   \cdot dW^{\nu}_{s}.
 \end{align*}
 Notice that $Z\in \Oc$,  $dt\x d\P$-a.e. and hence $  \delta(\nu)-\nu\cdot Z\ge 0,$ $dt\x d\P$-a.e. Then, it follows by comparison that 
\be \label{eq:minimal_sol}
	Y_{\sigma}\ge \Ec^{\nu}_{\sigma,T}[\xi],
	~~\mbox{for all}~ \nu \in \Uc ~\mbox{and}~~ \sigma \in \Tc.
\ee

\vspace{0.5em}
\noindent Conversely,  it is not difficult to deduce from the definition of $S$ that it satisfies a dynamic programming principle:
 \begin{equation*}
S(\sigma)=\esssup\left\{\Ec^{\nu}_{\sigma,\tau}[S(\tau)],\;\nu \in \Uc\right\},\;\forall\;( \sigma,\tau) \in \Tcd,
 \end{equation*}
 see e.g.~\cite{bouchard2014regularity}. 
 Taking $\nu\equiv 0$, we deduce that $S$ is a $\Ec^{0}$-supermartingale system. The existence of the aggregating process  $X \in \mathbf{X}^p_{\ell r}$  then follows from Theorem \ref{thm: doob-mertens decompo bsde}. Since it is also a  $\Ec^{\nu}$-supermartingale system for $\nu \in \Uc$, the same theorem implies that we can find $(  Z^{\nu},  A^{\nu})\in \Hb^{p}(\P^{\nu})\x\Ab^{p}(\P^{\nu})$ such that 
  \begin{equation*}
 X_{\sigma}
 = \xi+\int_{\sigma}^{T} ( g_{s}(X_{s},  Z^{\nu}_{s}) - \delta_{s}(\nu_{s}) )ds +  A^{\nu}_{T}-  A^{\nu}_{\sigma} -\int_{\sigma}^{T}   Z^{\nu}_{s} \cdot dW^{\nu}_{s},\; \;\sigma \in \Tc.
 \end{equation*}
 Identifying the quadratic variation terms implies that $  Z^{\nu}=  Z^{0}=:  Z$. Thus for all $\nu \in \Uc$,
 \begin{equation*}
 e(\nu):=\int_{0}^{T} ( \nu_{s}  Z_{s}- \delta_{s}(\nu_{s}) )ds \le \int_{0}^{T} ( \nu_{s}  Z_{s}- \delta(\nu_{s}) )ds +  A^{\nu}_{T}-  A^{\nu}_{0} =   A^{0}_{T}-  A^{0}_{0}.  
 \end{equation*}
We claim that if $N:=\{(\omega,t): Z_{t}(\omega)\notin \Oc_{t, \omega} \}$ has a non-zero measure w.r.t $d \P \x dt$, then we can find $\hat \nu \in \Uc$ such that $e(\hat \nu)\ge 0$ and $\P[e(\hat \nu)>0]>0$.
	However, for any real $\lambda>0$, one has $\lambda \hat \nu \in \Uc$ and  
	$e(\lambda \hat \nu)=\lambda e(\hat \nu)\le A^{0}_{T}-  A^{0}_{0}$, by the above, which is a contradiction since $A^0_T - A^0_0$ is independent of $\lambda$.
Hence, \reff{eq: constraint z} holds for $Z=Z^{0}$ and 
  \begin{equation*}
 X_{\sigma}
 = \xi+\int_{\sigma}^{T}  g_{s}(X_{s},  Z_{s})  ds +  A^{0}_{T}-  A^{0}_{\sigma} -\int_{\sigma}^{T}   Z_{s}\cdot dW_{s},\; \;\sigma \in \Tc.
 \end{equation*}
By \eqref{eq:minimal_sol}, it is clear that $(X, Z, A^0)$ is the minimal solution of \reff{eq: bsde supersol}-\reff{eq: constraint z}.
  
It remains to prove the above claim. Assume that $N$ has non-zero measure. Then, it follows from \cite[Thm.~13.1]{Rockafellar} that   $\{(\omega,t): \bar F_{t}(\omega):=\sup\{F_{t}(\omega,u),\;|u|=1\}\ge 2\iota\}$ has non-zero measure, for some $\iota>0$, in which 
$$F_{t}(\omega,u):=u\cdot Z_{t}(\omega)-\delta_{t}(\omega,u).$$
After possibly passing to another version (in the $dt\x d\P$-sense), we can assume that $Z$ is $\F^{\circ}$-progressively measurable. Since $\delta$ is $\Fc^{\circ}_{T}\otimes \Bc([0,T])\otimes \Bc(\R^{d})$-measurable,  $(\omega,t,u)\in \Omega\x [0,T]\x\R^{d} \longmapsto F_{t}(\omega,u) $ is Borel-measurable.  By \cite[Prop.~7.50 and Lem.~7.27]{BeSh78}, we can find a Borel map $(t,\omega)\longmapsto \hat u(t,\omega) $ such that  $|\hat u|=1$ and $F_{t}(\omega,\hat u(t,\omega))\ge \bar F_{t}(\omega)-\iota $ $dt\times d\P$-a.e.  Then, $\tilde u(t,\omega):=\hat u(t,\omega)\1_{N}(\omega,t)$ is Borel and satisfies $F_{t}(\omega,\hat u(t,\omega))\ge \iota\1_{N}(\omega,t)$  $dt\times d\P$-a.e. Since $F_{t}(\omega,\cdot)$ depends on $\omega$ only though $\omega_{\cdot\wedge t}$, the same holds for $(t,\omega)\longmapsto \hat u(t,\omega_{\cdot\wedge t})$, which is progressively measurable. We conclude by setting $$\hat \nu_{t}(\omega):=\hat u(t,\omega_{\cdot\wedge t})/(1+|\delta_{t}(\omega,\hat u(t,\omega_{\cdot\wedge t}))|).$$
 \ep

\appendix

\section{Appendix}

\subsection{Doob-Meyer decomposition for right-continuous supermartingales}
\label{subsec:RBSDE_DoobMeyer}

	We complete here the proof of Proposition \ref{prop:doob-general}, based on a personal communication with Nicole El Karoui.
	
	\vspace{0.5em}

\noindent {\bf Proof of Proposition \ref{prop:doob-general}.} 
 Let us start by considering the following reflected BSDEs with lower obstacle $X$ on $[0,\tau]$
\begin{equation}\label{eq:rbsde}
\begin{cases}
\displaystyle Y =Y_\tau+\int_{\cdot}^\tau g_s(Y_s,Z_s)ds-\int_\cdot^\tau Z_s\cdot dW_s-\int_\cdot^\tau dN_s-\int_\cdot^\tau dK_s,\\
\displaystyle Y\geq X\;\;   \mbox{on }  [0,\tau],\\
\displaystyle  \int_0^\tau(Y_{s-}-X_{s-})dK_s=0,
\end{cases}
\end{equation}
where $N$ is again a c\`adl\`ag martingale orthogonal to $W$, and $K$ is a c\`adl\`ag non-decreasing and predictable process. Since the obstacle $X$ is assumed to be c\`adl\`ag, the wellposedness of such an equation is guaranteed  by  \cite[Theorem 3.1]{BPTZ}.

\vspace{0.5em}
\noindent Let us now prove that we have $Y_t=X_t$, a.s., for any $t\in[0,\tau]$. Let us argue by contradiction and suppose that this equality does not hold. Without loss of generality, we can assume that $Y_0>X_0$ (otherwise, we replace $0$ by the first time when $Y>X+\iota$ for some $\iota>0$). Fix then some $\varepsilon>0$ and consider the following stopping time
$$\tau^\eps:=\inf\left\{t\geq 0,\ Y_t\leq X_t+\eps\right\}\wedge\tau.$$
Since $Y$ is strictly above $X$ before $\tau^\eps$, we know that $K$ is identically $0$ on $[\![0,\tau^\eps]\!]$, which implies that
$$Y_t=Y_{\tau^\eps}+\int_t^{\tau^\eps} g_s(Y_s,Z_s)ds-\int_t^{\tau^\eps} Z_s\cdot dW_s-\int_t^{\tau^\eps} dN_s.$$
Consider now the following BSDE on $[\![0,\tau^\eps]\!]$
$$y_t=X_{\tau^\eps}+\int_t^{\tau^\eps} g_s(y_s,z_s)ds-\int_t^{\tau^\eps} z_s\cdot dW_s-\int_t^{\tau^\eps} dn_s.$$
By standard a priori estimates (see for instance {\cite[Rem.~4.1]{BPTZ}}), we can find a constant $C>0$ independent of $\eps > 0$ s.t.
$$Y_0\leq y_0+C{\mathbb E\left[|X_{\tau^\eps}-Y_{\tau^\eps}|^{p}\right]^{\frac 1 p}}\leq y_0+C\eps.$$
But remember that $X$ is an $\mathcal E^{g}$-supermartingale, so that we must have $y_0\leq X_0$. Hence, we have obtained $Y_0\leq X_0+C\eps$, which implies a contradiction by arbitrariness of $\eps>0$.

\vspace{0.5em}
\noindent The uniqueness of the decomposition is then clear by identification of the local martingale part and the finite variation part of a semimartingale.
\ep

\subsection{Down-crossing lemma of $\Ec^g$-supermartingale} 
\label{subsec:downcroissing}

\noindent We provide here a down-crossing lemma for $\Ec^{g}$-supermartingales (defined in Section \ref{sec: doob meyer bsde} with $g$ satisfying \eqref{eq:LipCond_g} and \eqref{eq:integ} for some $p > 1$),
which is an extension of  Chen and Peng \cite[Thm 6]{chen2000general} 
(see also  Coquet et al.~\cite[Prop. 2.6]{coquet2002filtration}).
For completeness, we will also provide a proof. As in the classical case, $g\equiv 0$, it ensures the existence of right- and left-limits for $\Ec^{g}$-supermartingales, see Lemma \ref{lem: version limited a droite} below. 

\vspace{0.5em}
\noindent 
	For any $m>0$, we denote by $\Ec^{\pm m}_{\sigma, \tau}$ the non-linear expectation operator associated to the generator $(t,\omega, y ,z) \longmapsto \pm m |z|$ and stopping times $(\sigma, \tau) \in \Tcd$.
	Let $J := (\tau_n)_{n \in \N}$ be a countable family of stopping times taking values in $[0,T]$, 
	which are ordered, i.e. for any $i, j \in \N$, one has $\tau_i \le \tau_j$, {\rm a.s.},~or $\tau_i \ge \tau_j$, {\rm a.s.}
	Let  $a < b$, $X$ be some process and $J_n \subseteq J$ be a finite subset 
	($J_n = \{0 \le \tau_1 \le \cdots \le \tau_n \le T\}$).
	We denote by $D^b_a(X, J_n)$ the number of down-crossing of the process 
	$(X_{\tau_k})_{1 \le k \le n}$ from $b$ to $a$.
	We then define
	$$
		D_a^b(X, J)
		:=
		\sup \big\{
			D_a^b(X, J_n) ~: J_n \subseteq J, ~\mbox{and}~ J_n ~\mbox{is a finite set}
		\big\}.
	$$

\begin{Lemma}[Down-crossing]\label{lem: downcrossing} 
	Suppose that the generator $g$ satisfies \eqref{eq:LipCond_g}
	with Lipschitz constant $L$ in $y$ and $\mu$ in $z$, and \eqref{eq:integ} with $p > 1$.
	Let $X \in \mathbf{X}^p$ be a $\Ec^{g}$-supermartingale,  $J:= (\tau_n)_{n \in \N}$ be a countable family of stopping times taking values in $[0,T]$, which are ordered.
	Then, for all $a < b$,
	\begin{align} \label{eq:downcrossing}
		\Ec^{-\mu}_{0,T} \left[ D_a^b(X, J) \right]
		\le&\
		\frac{e^{LT}}{b-a} \Ec_{0,T}^{\mu} 
		\Big[ 
			e^{LT} (X_0 \wedge b -a)   - e^{-LT} ( X_T \wedge b -a)^+ \nonumber \\	
		&\hspace{4.2em}
			+ e^{LT} (X_T \wedge b -a)^- 
			+ e^{LT} \int_0^T | g_s(a,0)| ds
		\Big]. 
	\end{align}
\end{Lemma}
\proof
	First, without loss of generality, we can always suppose that $\tau^0 \equiv 0$ and $\tau^1 \equiv T$ belong to $J$,
	and also that $b > a = 0$.
	Indeed, whenever $b > a \neq 0$, 
	we can consider the barrier constants $(0, b - a)$, and the $\Ec^{\bar g}$-supermartingale $X- a$, with generator
	$
		\bar g_t(y, z) := g_t( y + a, z), 
	$
	which reduces the problem to the case $b > a = 0$.

\vspace{0.5em}
	\noindent Now, suppose that $J_n = \{ \tau_0, \tau_1, \cdots, \tau_n \}$ with $0 = \tau_0 < \tau_1 < \cdots < \tau_n =T$.
	We consider the following BSDE
	\begin{align*}
		y_t^i
		:&=
		X_{\tau_i}
		+
		\int_t^{\tau_i} g_s(y_s^i, z_s^i) ds 
		-
		\int_t^{\tau_i} z_s^i \cdot dW_s
		-
		\int_t^{\tau_i} dN^i_s \\
		&=
		X_{\tau_i}
		+
		\int_t^{\tau_i} \big( g_s(0,0) + \lambda^i_s y_s^i + \eta^i_s z_s^i  \big) ds 
		-
		\int_t^{\tau_i} z_s^i \cdot dW_s
		-
		\int_t^{\tau_i} dN^i_s,
	\end{align*}
	where $\lambda^i$ and $ \eta^i$ are  progressively measurable, coming from the linearization of $g$.
	In particular, we have $ |\lambda^i| \le L$ and $| \eta^i| \le  \mu$.
	Let us now consider another linear BSDE
	\begin{align} \label{eq:linearBSDE}
		\nonumber\bar y_t^i
		=&\
		X_{\tau_i}
		+
		\int_t^{\tau_i} \big(  - | g_s(0,0) | + \lambda^i_s \bar y_s^i +  \eta^i_s \bar z_s^i  \big) ds 
		-
		\int_t^{\tau_i} \bar z_s^i \cdot dW_s
		\\
		&-
		\int_t^{\tau_i} d \overline N^i_s.
	\end{align}
	By the comparison principle for BSDEs (see  \cite[Prop.~4]{kp}), and since $X$ is an $\Ec^{g}$-supermartingale,
	it is clear that 
	$$
		\bar y_{\tau_{i-1}}^i \le y_{\tau_{i-1}}^i \le X_{\tau_{i-1}}.
	$$
Solving the above linear BSDE \eqref{eq:linearBSDE}, it follows that
	\begin{align*}
		\bar y_{\tau_{i-1}}^i 
		&=
		\E^{\Q} \left[\left.
			X_{\tau_i} e^{\int_{\tau_{i-1}}^{\tau_i} \lambda^i_r dr} 
			-
			\int_{\tau_{i-1}}^{\tau_i} e^{\int_{\tau_{i-1}}^s \lambda^i_r dr} | g_s(0,0) | ds
		\right|
		\Fc_{\tau_{i-1}}
		\right],
	\end{align*}
	where $\Q$ is defined by
	\b*
		\frac{d\Q}{d \P}
		=
		 e^{-\frac12 \int_0^T  |\eta_s|^{2}   ds+
			\int_0^T  \eta_s\cdot  dW_s},
		&\mbox{with}&
		 \eta_s := \sum_{i=1}^n  \eta^i_s \1_{[\tau_{i-1}, \tau_i)}(s).
	\e*
	Let $\lambda_s := \sum_{i=1}^n \lambda^i_s \1_{[\tau_{i-1}, \tau_i)}(s)$,
	it follows that the discrete process $(Y_{\tau_i})_{0 \le i \le n}$ defined  by 
	\begin{equation*}
		Y_{\tau_i}
		:=
		X_{\tau_i} e^{\int_0^{\tau_i} \lambda_r dr} 
		-
		\int_0^{\tau_i} e^{\int_0^s \lambda_r dr} |g_s(0,0)  | ds
	\end{equation*}
	is a $\Q$-supermartingale.
	Define further
	\begin{equation*}
		\overline Y_{\tau_i}
		:=
		Y_{\tau_i}
		\wedge
		\left(b e^{LT} - \int_0^{\tau_i} e^{\int_0^s \lambda_r d\langle M \rangle_r} | g_s(0,0)  | ds \right),
	\end{equation*}
	which is clearly also a $\Q$-supermartingale.
	Let 
	$$
		u_t := b e^{\int_0^t \lambda_r dr} - \int_0^t e^{ \int_0^s \lambda_r dr} |g_s(0,0)  | ds,
	$$
	and
	$$
		l_t :=- \int_0^t e^{ \int_0^s \lambda_r dr} | g_s(0,0)  | ds.
	$$
	Denote then by $D_l^u(Y,J)$ (resp. $D_l^u(\overline Y, J)$) the number of down-crossing of the process $Y$ (resp. $\overline Y$) from the upper boundary $u$ to lower boundary $l$. It is clear that $D^u_l(Y,J) = D^u_l(\overline Y, J)$.
	Notice that $l_t$ is decreasing in $t$,
	so that we can apply the classical down-crossing theorem for supermartingales (see e.g.~Doob \cite[p.446]{Doob83}) to $\overline Y$,
	and obtain that
	\begin{align*}
		&\E^{\Q} \left[ D^b_0(X, J)\right]\\
		&\le
		\E^{\Q} \left[ D^u_l (\overline Y, J) \right]\\
		&\le
		\frac{e^{LT} }{b} \E^{\Q} \left[ (\overline Y_0 - \overline Y_T ) - (u_T - \overline Y_T) \wedge 0 \right] \\
		&\le
		\frac{e^{LT} }{b}  \E^{\Q} 
		\left[
			X_0 \wedge (b e^{LT}) -e^{\int_0^T \lambda_s ds} (X_T \wedge b)+~e^{LT} \int_0^T   | g_s(0,0) | ds
		\right].
	\end{align*}
	Notice that $|\lambda_s| \le L$, $| \eta_s| \le \mu$ and $(X_T \wedge b) = (X_T \wedge b)^+ - (X_T \wedge b)^-$. Therefore, we have proved \eqref{eq:downcrossing} for the case $b > a = 0$, from which we conclude the proof, by our earlier discussion.
\ep

\begin{Lemma}\label{lem: version limited a droite} Let $X \in \mathbf{X}^p$ be a $\Ec^{g}$-supermartingale of class $(D)$. Then, it admits right- and left-limits outside an evanescent set. 
\end{Lemma}

\proof We follow well-known arguments for (classical) supermartingales. Let $(\vartheta_{n})_{n}\subset \Tc$ be a non-increasing sequence of stopping times. Then, $(X_{\vartheta_{n}})_{n\ge 1}$ converges $\as$ This is an immediate consequence of the down-crossing inequality of Lemma \ref{lem: downcrossing}, see e.g.~\cite[Proof of Thm V-28]{DellacherieMeyer}. Set $\bar X:=X/(1+|X|)$. Then, \cite[Thm VI-48]{DellacherieMeyer} implies that, up to an evanescent set, $\bar X$ admits right-limits. Since $a/(1+|a|)=b/(1+|b|)$ implies $a=b$, for all $a,b\in \R$, this shows that $ X$ admits right-limits, up to an evanescent set. The existence of left-limits is proved similarly by considering non-decreasing sequences of stopping times.\ep

\bibliographystyle{plain}

  \end{document}